\documentclass[11pt,a4paper]{article}
\usepackage[T1]{fontenc}
\usepackage{lmodern}
\usepackage[margin=1in]{geometry}
\usepackage{amsmath,amssymb,amsthm,booktabs}
\usepackage[hidelinks]{hyperref}
\hypersetup{pdftitle={Two Knuth orbits of order 256 without primitive semifield coordinatizations},
  pdfauthor={I. F. Rúa and E. F. Combarro}}
\newtheorem*{result}{Theorem}
\newcommand{\F}{\mathbb F}
\newcommand{\Aut}{\operatorname{Aut}}
\newcommand{\Atp}{\operatorname{Atp}}

\title{Two Knuth orbits of order 256 without\\ primitive semifield coordinatizations}
\author{I. F. R\'ua\thanks{Department of Mathematics, University of Oviedo, Asturias, Spain.
  Email: \texttt{rua@uniovi.es}.}
\and E. F. Combarro\thanks{Department of Computer Science, University of Oviedo, Asturias, Spain.}}
\date{}

\begin{document}
\maketitle
\begin{abstract}
We exhibit two Knuth orbits of semifield planes of order \(256\) in which
every coordinatizing semifield is nonprimitive on both the left and the right.
Each orbit consists of two isotopy classes. The first orbit contains \(168\)
isomorphism classes of unital semifields and the second contains \(32\).
We give binary multiplication matrices for a representative of each orbit
and an exhaustive computation that verifies nonprimitivity and the counts
directly from these matrices.
\end{abstract}

\section{The two orbits}

In a finite semifield \(S\) with unity \(1\), write \(R_a(x)=xa\) and
\(L_a(x)=ax\). An element \(a\) is \emph{right primitive} if its right
principal powers \(1,a,aa,(aa)a,\ldots\) exhaust \(S^\times=S\setminus\{0\}\).
Left primitivity is defined using the iterates of \(L_a\) at \(1\).
A semifield is \emph{two-sided nonprimitive} if it has neither a right nor
a left primitive element.

Wene studied primitivity in finite semifields~\cite{Wene}. The commutative
binary Knuth semifield of order \(32\) is two-sided nonprimitive~\cite{Rua},
as is an example of order \(64\) found by Hentzel and R\'ua~\cite{HR}.
Both examples are discussed by Kravtsova and Kuzmin~\cite{KK}.
Primitivity can change under isotopy~\cite{Rua,HR}. To rule out a primitive
coordinatization for a plane, one must therefore consider every unital
semifield in its isotopy class~\cite{Albert}.
The six permutations of the indices of a multiplication tensor give the
associated Knuth orbit of isotopy classes~\cite{Knuth}.

\begin{result}
There are two distinct Knuth orbits of order \(256\), each consisting of
two isotopy classes, such that none of the four corresponding planes
admits a left or a right primitive semifield coordinatization.
Together they contain \(200\) pairwise nonisomorphic unital semifields.
\end{result}

Table~\ref{tab:invariants} gives the characteristics of the representatives
\(S_1,S_2\) specified below. The action in \(C_3\rtimes C_4\) is by
inversion, and \(Q_8\) denotes the quaternion group of order eight.
The distinct nucleus orders separate the two Knuth orbits, since isotopy
preserves the orders of the nuclei and tensor permutations permute them.

\begin{table}[ht]
\centering\small
\caption{Invariants of the representatives and numbers of unital isomorphism classes.}
\label{tab:invariants}
\begin{tabular}{@{}lccccc@{}}
\toprule
& \((|N_l|,|N_m|,|N_r|)\) & \(\Aut(S_i)\) & \(|\Atp(S_i)|\) &
\multicolumn{2}{c}{Isomorphism classes}\\
\cmidrule(l){5-6}
& & & & per plane & per Knuth orbit\\
\midrule
\(S_1\) & \((2,2,2)\) & \(C_3\rtimes C_4\) & \(900\)  & \(84\) & \(168\)\\
\(S_2\) & \((4,4,4)\) & \(Q_8\)            & \(5400\) & \(16\) & \(32\)\\
\bottomrule
\end{tabular}
\end{table}

\section{Multiplication matrices}

Let \(V=\F_2^8\) with ordered basis \(e_0,\ldots,e_7\). The label of
\(\sum_i x_i e_i\) is \(\sum_i 2^i x_i\), written as a two-digit
hexadecimal number. For each representative, row \(j\) of
Table~\ref{tab:matrices} lists the eight \emph{columns} of \(A_j=R_{e_j}\)
in order. Thus the \(i\)-th entry, counting from zero, is \(e_i*e_j\).
For column vectors \(x\) and \(y=\sum_j y_je_j\) the product is
\[
                 x*y=\Bigl(\sum_{j=0}^7 y_j A_j\Bigr)x .
\]
In each example \(A_0=I_8\) and \(A_j e_0=e_j\), so the unity is \(e_0\),
labelled \texttt{01}. All \(255\) nonzero right multiplication matrices
are invertible, as can be checked directly from the table. This proves the
semifield property.

\begin{table}[ht]
\centering\small
\caption{Column labels of the binary right multiplication matrices \(A_j\).}
\label{tab:matrices}
\begin{tabular}{@{}c@{\qquad}l@{\qquad}l@{}}
\toprule
\(j\) & \(S_1\) & \(S_2\)\\
\midrule
0 & \texttt{01 02 04 08 10 20 40 80} & \texttt{01 02 04 08 10 20 40 80} \\
1 & \texttt{02 04 08 03 50 A0 70 E0} & \texttt{02 04 08 10 20 40 80 81} \\
2 & \texttt{04 08 03 06 20 40 80 30} & \texttt{04 33 C4 F5 98 9F 1B 3D} \\
3 & \texttt{08 03 06 0C A0 70 E0 F0} & \texttt{08 3E F4 65 B0 AD DA DC} \\
4 & \texttt{10 20 40 80 11 43 35 CF} & \texttt{10 93 7B B4 85 47 39 46} \\
5 & \texttt{20 40 80 30 55 7F F2 96} & \texttt{20 C7 03 7E 34 71 7F 1E} \\
6 & \texttt{40 80 30 60 22 86 6A BD} & \texttt{40 0C DA B0 0E 10 85 B6} \\
7 & \texttt{80 30 60 C0 AA ED D4 1C} & \texttt{80 EB 5C E3 4E 6E 48 93} \\
\bottomrule
\end{tabular}
\end{table}

There is also a field description of the first table. Let
\(K=\F_2[\eta]\), with \(\eta^4+\eta+1=0\). In the basis
\((1,\eta,\eta^2,\eta^3)\) of each summand of \(K\oplus K\), the product is
\[
 (u,w)*(v,z)=(uv+w^4z^8,\;uz+wv^8+w^2z^8).
\]
It belongs to Knuth's quadratic construction~\cite{Knuth}.

\section{Finite verification}

For any of the six tensor permutations, denote its product on \(V\) by
\(*\) and its multiplication maps by \(R_y,L_x\). For \(a,b\ne0\)
the principal isotope
\[
                 x\circ y=R_b^{-1}(x)*L_a^{-1}(y)
\]
has unity \(a*b\). Every unital isotope is isomorphic to one of this form.
As the multiplier varies over \(V^\times\), its right and left multiplication
maps run through the respective families
\begin{equation}\label{eq:quotients}
        R_yR_b^{-1}\quad(y\ne0),\qquad L_xL_a^{-1}\quad(x\ne0).
\end{equation}
Indeed, right multiplication by \(u\) in the isotope is
\(R_{L_a^{-1}(u)}R_b^{-1}\). The left-hand analogue is
\(L_{R_b^{-1}(u)}L_a^{-1}\).

A permutation of the \(255\) nonzero vectors has a cycle of length \(255\)
if and only if the orbit of \(e_0\) has this length. Thus one may start
at \texttt{01} for every test, including isotopes whose unity differs
from \(e_0\). For each side, the computation forms every ordered pair of
nonzero parameters in~\eqref{eq:quotients} and iterates the corresponding
map until the starting vector returns. For each of the six tensor
permutations of each example, the number of full cycles is zero in
all \(255^2=65\,025\) tests on the right and in all \(65\,025\) tests on
the left.

To determine the Knuth orbits we enumerate isotopies from \(S_i\) to
each tensor permutation \(P\). For \(g=\texttt{12}\) in \(S_1\) and
\(g=\texttt{02}\) in \(S_2\), the vectors \(R_g^k(1)\), \(0\leq k<8\),
form a binary basis. An isotopy \((F,G,H)\) satisfies
\(H(x*y)=F(x)\mathbin{\cdot}G(y)\), where \(\cdot\) is the product of \(P\).
Write \(a=F(1), b=G(1), y=G(g)\). With multiplication maps in \(P\)
distinguished by a superscript, we have
\[
 H(1)=a\cdot b,\qquad HR_g=R_y^P(R_b^P)^{-1}H,\qquad
 F=(R_b^P)^{-1}H,\quad G=(L_a^P)^{-1}H.
\]
Each choice of the three nonzero images therefore determines a unique
candidate on the power basis. Invertibility and the \(64\) basis products
decide whether it is an isotopy. For \(S_1\) the count is \(900\) for
each even permutation and zero for each odd permutation. The corresponding
counts for \(S_2\) are \(5400\) and zero. Hence in each orbit the two
classes are represented by the original multiplication and its opposite.

Finally \(\Atp(S_i)\) acts on the pairs defining principal isotopes by
\((a,b)\mapsto(F(a),G(b))\). Its orbits are the isomorphism classes of
these isotopes. Indeed, \(H\) identifies the two corresponding products.
Conversely, an isomorphism \(H\) from the isotope indexed by \((a,b)\) to
that indexed by \((c,d)\) gives the autotopy
\((R_d^{-1}HR_b,L_c^{-1}HL_a,H)\).
The sums of the fixed-pair counts in Burnside's lemma are \(75\,600\)
and \(86\,400\), giving \(75\,600/900=84\) and \(86\,400/5400=16\).
Reversing multiplication gives the same counts for the opposite classes.
This yields \(2(84+16)=200\).

We used the multiplication matrices in Table~\ref{tab:matrices} to compute
the invariants in Table~\ref{tab:invariants} and to verify that none of the
maps in~\eqref{eq:quotients} has a cycle of length~\(255\).
The C++17 source code (\texttt{verify.cpp}) and its output (\texttt{results.txt})
are included as supplementary material.

\section{Final remarks}

This note is a preliminary announcement. In an extended version we will
use quartic field norms to explain why every semifield in either Knuth
orbit is nonprimitive on both sides.

\paragraph{Acknowledgments.}
This work was partially supported by grant PID2023-146520OB-C22 funded by
MCIU/AEI/10.13039/501100011033/FEDER, EU, and by grant
IDE/2024/000734 funded by the Principality of Asturias. Further support
came from the Ministry for Digital Transformation and Civil Service of
the Spanish Government through the QUANTUMENIA project under the Quantum
Spain programme, and from the European Union through the Recovery,
Transformation and Resilience Plan---NextGenerationEU, within the
framework of the Digital Spain 2026 Agenda. The first author is a member
of the Universidad de Oviedo research team GACYC and has also been
supported by the Spanish Network of Mathematics in the Information
Society (MatSI). We used ChatGPT (OpenAI) to help draft this note and write
the C++17 verification program.


\begin{thebibliography}{9}
\bibitem{Albert}
A. A. Albert, Finite division algebras and finite planes,
in \emph{Combinatorial Analysis}, Proc. Sympos. Appl. Math. \textbf{10},
Amer. Math. Soc. (1960), 53--70.
\bibitem{HR}
I. R. Hentzel and I. F. R\'ua,
\href{https://doi.org/10.1142/S0218196707004220}{Primitivity of finite semifields with 64 and 81 elements},
\emph{Internat. J. Algebra Comput.} \textbf{17} (2007), 1411--1429.
\bibitem{Knuth}
D. E. Knuth,
\href{https://doi.org/10.1016/0021-8693(65)90018-9}{Finite semifields and projective planes},
\emph{J. Algebra} \textbf{2} (1965), 182--217.
\bibitem{KK}
O. V. Kravtsova and I. K. Kuzmin,
\href{https://www.mathnet.ru/eng/jsfu1220}{On spectra and minimal polynomials in finite semifields},
\emph{J. Sib. Fed. Univ. Math. Phys.} \textbf{18} (2025), 41--50.
\bibitem{Rua}
I. F. R\'ua,
\href{https://doi.org/10.1081/AGB-120027930}{Primitive and non primitive finite semifields},
\emph{Comm. Algebra} \textbf{32} (2004), 793--803.
\bibitem{Wene}
G. P. Wene,
\href{https://doi.org/10.1007/BF02227457}{On the multiplicative structure of finite division rings},
\emph{Aequationes Math.} \textbf{41} (1991), 222--233.
\end{thebibliography}
\end{document}